# A convex model for induction motor starting transients imbedded in an OPF-based optimization problem

H. Sekhavatmanesh, *Student Member, IEEE*
R. Cherkaoui, *Senior Member, IEEE*
Power System Research Group
École Polytechnique Fédérale de Lausanne EPFL
Lausanne, Switzerland
hossein.sekhavat@epfl.ch

J. Rodrigues[1]
C. L. Moreira[1,2]
J.A.P. Lopes[1,2], *Fellow Member, IEEE*
[1]INESC Technology and Science
[2]Faculty of Engineering of the University of Porto
Porto, Portugal

*Abstract*— Large horsepower induction motors play a critical role in the operation of industrial facilities. In this respect, the distribution network operators dedicate a high priority to the operational safety of these motor loads. In this paper, the induction motor starting is modeled analytically and in a semi-static fashion. This model is imbedded in a convex distribution system restoration problem. In this optimization problem, it is aimed to determine the optimal status of static loads and the optimal dispatch of distributed generators such that: a) the induction motors can be reaccelerated in a safe way and, b) the total power of static loads that cannot be supplied before the motor energization, is minimized. The proposed optimization problem is applied in the case of a distribution network under different simulation scenarios. The feasibility and accuracy of the obtained results are validated using a) off-line time-domain simulations, and b) Power Hardware-In-the-Loop experiments.

## I. INTRODUCTION

Large horsepower induction motors play critical roles in most of industrial facilities [1]. Therefore, re-energization of these loads, following a fault, has a huge priority for Distribution Network Operators (DNOs). When a fault occurs in a distribution network, once it is isolated, the unsupplied area (named as *off-outage area*) is transferred to the healthy neighboring feeders through switching operations. Under this new configuration, the loads in the off-outage area should be re-energized in a certain sequence such that the total energy not supplied will be minimized. In this regard, the re-energization of induction motors should be done with special precaution in order to avoid violation of operational constraints I) in the motor side and II) in the network side.

In the motor side, the large inrush current of the induction motor could make extreme voltage drop at the motor connecting node and ultimately a reduction in the electromagnetic torque. This may result in either the motor stall or a prolonged acceleration that could trigger the tripping of the overheating protection at the motor side. Regarding from the network side, the starting of the induction motor may cause the tripping of over-current or under-voltage protections of healthy lines and nodes in the distribution network. In order to avoid these risks, a reliable restoration plan should be derived while considering the dynamic model of the induction motor and the detailed model of the distribution network.

The analysis of motor acceleration transients and its application to different levels of power system operation has been the subject of many papers in the literature. One group of papers study the starting of the motor loads following a complete-blackout in the power system. In this regard, the starting transients of the large induction motors are evaluated during different stages of the power system restoration [3]-[5]. Regarding the low-voltage networks, the authors in [6]–[8] evaluate the operational safety of an islanded microgrid in terms of the frequency and voltage stability in case of induction motor starting. The re-acceleration scheme of motor drives in a small industrial facility following an outage is studied in [1], [9]. Apart from the permanent fault situations, there are papers evaluating the behavior of induction motors in other situations, such as voltage sags [10] and short interruptions [11].

In all the papers mentioned above, the behavior of the induction motor is assessed using time-domain simulations only. However, if these dynamic simulations are used for decision-making purposes, different combination of decision variables (e.g. load energization sequences) should be tested separately. This could be very time-consuming for the problems with a huge and complex solution space. This is the case for the problem of distribution network restoration studied in this paper. In fact, it needs to formulate analytically the motor starting dynamic and integrate it into an optimization problem.

Such analytical approaches are proposed in [12]–[15]. The authors of [12] developed nonlinear differential-algebraic formulations to estimate the voltage dip during the motor starting. This voltage dip is predicted in [13] using neural network for an induction motor with a certain kVA capacity installed on a bus with a certain short circuit capacity. The developed formulations are all nonlinear and non-convex as well. Therefore, they cannot be integrated into the convex optimization problems. The presented formulation in [14] is incorporated into a maximum restorable load problem and solved iteratively using a heuristic approach. This approach is applicable only in the case of simple problems involving only one decision variable. The authors in [15] presented a quadratic optimization problem for the minimization of voltage deviation with respect to the nominal value during the motor starting. In that paper, a model predictive control approach is used to co-ordinate all the reactive power resources in the system in order to support the voltage profile during the motor acceleration. According to this method, the reactive power demanded by the motor loads during the acceleration period is approximated as a simple algebraic function of terminal voltage magnitude. This simplified model does not account for the dynamics of the motor starting. In consequence, it cannot represent correctly the operational safety constraints.

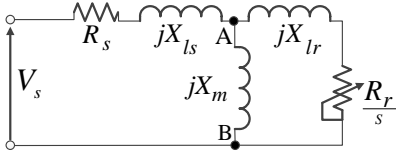

Fig. 1. The equivalent circuit of the induction motor.

In this paper, an analytical optimization model is derived for the starting of induction motors in a *semi-static* fashion. This formulation is integrated into a distribution network restoration problem resulting in a convex mixed-integer optimization problem. In this restoration problem, it is aimed to determine I) the optimal status of the static loads and II) the optimal dispatch of DGs such that the induction motor can be re-accelerated in a safe way. More precisely, it is intended to minimize the total power of static loads -considering their priorities- that cannot be supplied before the motor energization. The proposed optimization problem is tested thanks to a distribution network under different simulation scenarios. The feasibility of the obtained solution is verified using I) off-line time-domain simulations and II) using a Power-Hardware-In-the-Loop (PHIL) test experiment at INESC TEC.

## II. SEMI-STATIC MODEL OF THE INDUCTION MOTOR STARTING

In this section, an analytical and semi-static model is proposed for the induction motor starting such that it can be integrated into a convex OPF problem. Fig. 1 shows the equivalent circuit of the induction motor, neglecting the dynamics of rotor fluxes. Using this equivalent circuit, the electrical torque of the induction motor ($T^{ele}$) is formulated in the following as a function of the rotor slip:

$$T^{ele} = \frac{R_r V_{th}/S}{(R_r/S + R_{th})^2 + (X_{lr} + X_{th})^2} \quad (1)$$

where:

$V_{th} = \frac{V_s X_m^2}{R_s^2 + (X_{ls}^2 + X_m^2)}$

$R_{th} = \frac{R_s X_m^2}{R_s^2 + (X_{ls}^2 + X_m^2)}$

$X_{th} = \frac{R_s^2 X_m + X_{ls} X_m (X_{ls} + X_m)}{R_s^2 + (X_{ls}^2 + X_m^2)}$

$V_s$ is the square of the stator voltage. S is the motor slip. $R_s$ and $X_{ls}$ represent the resistance and reactance of the stator. $R_r$ and $X_{lr}$ are the resistance and reactance of the rotor. $X_m$ is the magnetization reactance. $V_{th}$, $R_{th}$ and $X_{th}$ are the Thevenin's voltage square, resistance and reactance seen from the rotor terminals (AB in Fig. 1), respectively.

Fig. 2 shows the typical electrical torque-slip curve of the induction motor assuming that the motor terminal voltage is fixed. The load torque ($T^{load}$) refers to the summation of the mechanical torque on the shaft ($T^{mec}$) and the friction and windage toque ($K_D(1 - S_k)$). In Fig. 2, it is assumed that the mechanical load has a linear toque-speed characteristic. The different between the electrical and load torques is named the acceleration torque and referred as $T^{acc}$ in Fig. 2.

When starting, the slip is unity and then it decreases gradually to a stable point where it is close to 0 and the speed is close to the synchronous speed. As it can be seen, the electromagnetic torque is a non-linear function of the motor slip. In order to build a convex model of the induction motor starting, we divide the slip range between standstill ($s = 1$) and the stable point into $K_{max}$ fixed steps with equal length of $\Delta s$. From hereafter, slip steps (or simply steps) refer to these short and sequential intervals of the motor acceleration period. During each step $k$, the slip is assumed to be fixed and equal to a given value ($S_k$). The justification for this assumption is provided in Appendix A.

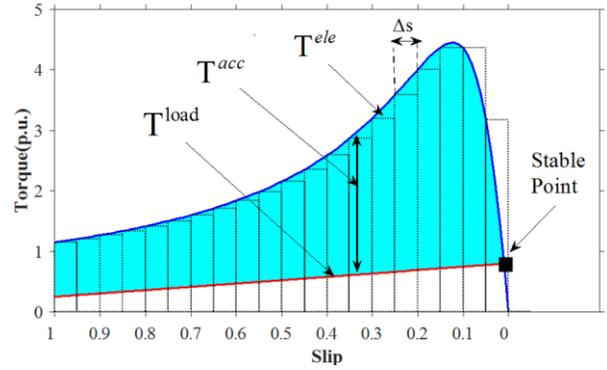

Fig. 2. The discretized electrical, mechanical, and accelerating torques of the induction motor

Given the slip value $S_k$ at each step $k$ as a parameter, the electrical torque of the induction motor can be obtained using (1) as a linear function of the square of voltage terminal ($V_s$), which is a state variable in the OPF problem. Now, we want to derive the time length of a given step $k$ indicated by $\Delta t_k$. This is defined as the period during which the slip changes from $S_k$ to $S_{k+1} = S_k - \Delta S$. For a given step $k$, $\Delta t_k$ is derived from (2) after some manipulations on the dynamic motion equation given in (20). It should be noted that since the acceleration torque is assumed to be fixed during each step, the time derivative of the slip can be approximated as $\frac{\Delta S}{\Delta t}$.

$$\frac{1}{\Delta t_k} = \frac{1}{2H \cdot \Delta S}\left(T_k^{ele} - T_k^{mec} - K_D(1 - S_k)\right) \quad (2)$$

As formulated in (2), the inverse of $\Delta t_k$ is a linear function of the accelerating torque at step k. The time length of each single step k ($\Delta t_k$) is obtained using a piece-wise linear approximation of an inverse function in the form of $f(x) = \frac{1}{x}$. The formulation of this piecewise linear approximation is explained in Appendix III.A for an arbitrary continuous function. For a given step k, the acceleration time is defined as the time that is taken for the slip to change from standstill ($s = 1$) to the value at step k ($s = S_k$). The acceleration time of step k is indicated by $t_k$ and obtained by summing up the time lengths of all previous steps as expressed in (3). This acceleration time will be used in section III in order to derive the transient voltage and current limits.

$$t_k = \sum_{k^*=1}^{k} \Delta t_{k^*} \quad (3)$$

While fixing the slip during each step, all the parameters of the motor equivalent circuit shown in Fig. 1 will be fixed too. Therefore, the network can be represented in phasor domain for each single step. It should be mentioned that the DC term imbedded in the starting current of the induction motor is neglected in this paper. This assumption is justified because of the low X/R ratio in distribution networks. Considering this resistive characteristic, the DC term of the starting current disappears during the first instants of the acceleration period. With this assumption, we incorporate the AC power flow equations into the above-mentioned discretization method. The aim is to obtain the electrical state variables in the whole distribution network for each step.

### III. RESTORATION PROBLEM FORMULATION

In this section, the load restoration problem is presented as an example to show how the proposed semi-static model of the motor starting could be included into a convex optimization problem. It is assumed that the fault is already isolated and the backbone of the off-outage area is energized by healthy neighboring feeders. According to the restoration strategy, these actions are followed by the re-energization of the unsupplied loads in an optimal sequence. In this paper, it is assumed that there is only one induction motor load in the off-outage area. In addition, it is assumed that the energization time of this single motor load is pre-defined. In the proposed optimization problem, the main decision variables are two-folds, namely I) binary variables $L_i$ which indicate if the static load at node $i$ is energized or not and II) continuous variables $P_i^{DG}/Q_i^{DG}$ which represent the active-reactive set points for the DG at node $i$ during the motor acceleration. The obtained results of the optimization problem are deployed in the network before energizing the induction motor. The following presents the formulation of the optimization problem in the form of a Mixed-Integer Second-Order Cone Programming (MISOCP).

Minimize:    $F^{obj} = W_{re}.F^{re} + W_{op}.F^{op}$    (4)

$$F^{re} = \sum_{i \in N} D_i.(1 - L_i).P_i^0 \tag{5}$$

$$F^{Op} = \sum_{k=1}^{k_{max}} \sum_{ij \in W} r_{ij}.F_{ij,k} \tag{6}$$

Subject to:

$$L_i = 1 \qquad i \in N \backslash N^* \tag{7}$$

$$P_{i,k}^D = \begin{cases} P_{i,k}^0\left(1 + (V_{i,k} - 1)\right)^{kp_i/2} \approx P_{i,k}^0(1 + \frac{kp_i}{2}(V_{i,k} - 1) : i \in N_s \\ V_{i,k}.\dfrac{R_{T,k}}{\sqrt{R_{T,k}^2 + X_{T,k}^2}} \qquad : i \in N_m^* \\ P_i^0 \qquad : i \in N_m \backslash N_m^* \end{cases} \tag{8}$$

$$Q_{i,k}^D = \begin{cases} Q_{i,k}^0\left(1 + (V_{i,k} - 1)\right)^{kq_i/2} \approx Q_{i,k}^0(1 + \frac{kq_i}{2}(V_{i,k} - 1) : i \in N_s \\ V_{i,k}.\dfrac{X_{T,k}}{\sqrt{R_{T,k}^2 + X_{T,k}^2}} \qquad : i \in N_m^* \\ Q_i^0 \qquad : i \in N_m \backslash N_m^* \end{cases} \tag{9}$$

$$V_{i,k} - V_{j,k} - 2(r_{ij}.p_{ij,k} + x_{ij}.q_{ij,k}) = 0 \quad \forall ij \in W \tag{10}$$

$$p_{ij,k} = (\sum_{\substack{i^* \neq i \\ (i^*,j) \in W}} p_{ji^*,k}) + r_{ij}.F_{ij} + L_j.P_{j,k}^D - (P_{j,k}^{Sub} + P_{j,k}^{DG}) \quad \forall ij \in W \tag{11}$$

$$Q_{ij,k} = (\sum_{\substack{i^* \neq i \\ (i^*,j) \in W}} Q_{ji^*,k}) + x_{ij}.F_{ij} + L_j.Q_{j,k}^D - (Q_{j,k}^{Sub} + Q_{j,k}^{DG}) \quad \forall ij \in W \tag{12}$$

$$\left\| \begin{matrix} 2p_{ij,k} \\ 2q_{ij,k} \\ F_{ij,k} - V_{i,k} \end{matrix} \right\|_2 \leq F_{ij,k} + V_{i,k} \quad \forall ij \in W \tag{13}$$

$$0 \leq P_{i,k}^{DG} \leq P_{i,max}^{DG}$$
$$Q_{i,min}^{DG} \leq Q_{i,k}^{DG} \leq Q_{i,max}^{DG} \qquad \forall i \in N_{DG} \tag{14}$$

$$\left\| \begin{matrix} P_{i,k}^{DG} \\ Q_{i,k}^{DG} \end{matrix} \right\|_2 \leq S_{i,max}^{DG} \qquad \forall i \in N_{DG} \tag{15}$$

$$T_k^{ele} \geq T_k^{mec} + K_D(1 - S_k) \tag{16}$$

$$L_i.\tilde{V}_k^{min} \leq V_{i,k} \qquad \forall i \in N_p \tag{17}$$

$$F_{ij,k} \leq \tilde{F}_{ij,k}^{max} \qquad \forall ij \in L_p \tag{18}$$

As formulated in (4), the objective function ($F^{obj}$) includes reliability ($F^{re}$) and operational ($F^{op}$) objective terms. The main objective of the restoration problem is the reliability term. As formulated in (5), this term tends to minimize the total costs associated with the loads that cannot be supplied ($L_i = 0$) before the induction motor energization. The cost of not supplying a given load at node $i$ is represented by its nominal active power ($P_i^0$) multiplied by its priority factor ($D_i$). Unlike the reliability term, the operational term has a very small weighting coefficient. This term is expressed in (6) as the total active power losses in the distribution network, where $r_{ij}$ and $F_{ij,k}$ represent the resistance and the square of current magnitude of line $ij$, respectively. This term is included in the objective function just to ensure the exactness of the relaxed AC-OPF formulation. In order to deploy the predefined hierarchical priority between the two optimization criteria, the weighting factors $W_{re}$ and $W_{op}$ are applied on the normalized terms of expressions formulated in (5) and (6), respectively. The index k represents the predefined steps. Each constraint with the index of k should be hold for every step $k \in \{1,2,..,k_{max}\}$.

Constraint (7) indicates that loads that are not in the off-outage area ($N^*$) should be remained supplied. The active ($P_{i,k}^D$) and reactive ($Q_{i,k}^D$) power of a given load at node i and step k are formulated in (8) and (9), respectively. The set of static loads is represented by $N_s$. These loads are formulated using the exponential model, where $K_p$ and $K_q$ represent the corresponding active and reactive coefficients, respectively. This model is linearized using the binomial approximation, assuming that $V_{i,k}$ is close to $1\ p.u.$ The set of all the motor loads is represented by $N_m$. It is assumed that the motor loads that are already energized consume their nominal active ($P_i^0$) and reactive ($Q_i^0$) powers. $N_m^*$ represents the index of the induction motor that is unsupplied and is subject to the re-energization. The active and reactive powers of this motor load are derived according to the equivalent circuit shown in Fig. 1. $R_{T,k}$ and $X_{T,k}$ represent the Thevenin's equivalent resistance and reactance seen from the motor terminals at a given step k. As it can be seen, the active and reactive powers of all the loads at a given step $k$ are represented as linear functions of the terminal voltage square ($V_{i,k}$), which is an electrical state variable.

Constraints (10)-(15) represent the relaxed formulation of the AC-OPF problem [2]. The aim is to obtain the electrical state variables of the whole distribution network for every step k during the motor acceleration period. Constraint (10) expresses the nodal voltage equation, where $p_{ij,k}/q_{ij,k}$ are the active/reactive power flow on line ij and at step k and $r_{ij}/x_{ij}$ are the resistance/ reactance of line ij. Constraints (11) and (12), respectively, concern with the active and reactive power balances at the extremities of each line, where $P_{i,k}^{Sub}/Q_{i,k}^{Sub}$ are the active/reactive power injection at the substation node i (slack bus) and at step k. Constraint (13) is the relaxed version of the current flow equation in each line according to [2]. Constraint (14) limits the active and reactive power injections of each DG, while its apparent power is limited by the cone constraint given in (15). $P_{i,max}^{DG}$, $Q_{i,max}^{DG}$, and $S_{i,max}^{DG}$ are the active, reactive and apparent power capacity of the DG at node *i*, respectively. The set of nodes which host a DG are represented by $N_{DG}$. The products of the binary variable $L_i$ and the positive continuous variables $P_{i,k}^D$ and $Q_{i,k}^D$ introduce non-linear terms in (11) and (12). In order to preserve the linearity, these terms are re-formulated according to the strategy given in Appendix C.

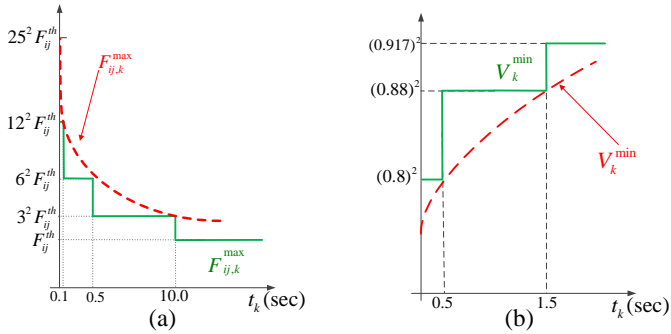

Fig. 3. Typical protection curve of a) over-current and b) under-voltage relays used in this paper.

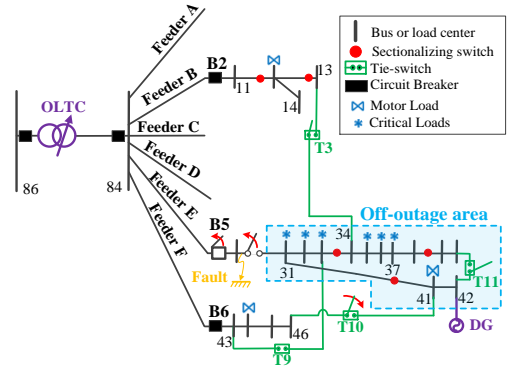

Fig. 4. Part of the test distribution network under the post-fault configuration.

Constraints (16)-(18) express the transient constraints regarding the safe starting of the induction motor in a distribution network. In order to avoid the induction motor to stall, (16) enforces the electrical torque to be larger than or equal to the load torque. For a given slip $S_k$, the electrical torque is obtained using (1) and the mechanical load torque is determined according to the torque-speed curve of the mechanical load. This curve is assumed to be given for a specific load on the shaft. Constraints (17) and (18) represent, for every step k, the function of under-voltage and over-current protection devices protecting the loads at nodes $N_p$ and protecting the lines $L_p$, respectively. Fig. 3 shows the typical protection curves reported in [3] and [4]. As it can be seen, the values of the under-voltage ($V_k^{min}$) and over-current ($F_{ij,k}^{max}$) limits for a given step k are functions of the acceleration time $t_k$ at step k, that is already formulated in (3). In Fig. 3.a, $F_{ij}^{th}$ represents the square of the nominal thermal ampacity limit of line ij. In order to preserve the linearity in terms of variable $t_k$, $V_k^{min}$ and $F_{ij,k}^{max}$ are approximated by $\tilde{V}_k^{min}$ and $\tilde{F}_{ij,k}^{max}$, respectively, according to the piecewise linear approximation method that is explained in Appendix B. In order to preserve the linearity of (17), the product of binary variable $L_i$ and positive variable $\tilde{V}_k^{min}$ is re-formulated according to the linearization strategy explained in Appendix C.

## IV. SIMULATION RESULTS

In order to validate the functionality of the proposed restoration problem, it is tested on the model of a 11.4 kV distribution network in Taiwan [18]. Fig. 4 shows only a part of this network that is affected by a fault on line 30-31. Once the fault is isolated by opening breaker B5 and switch 30-31, the rest of feeder E (off-outage area) is restored by closing tie-switch T10. As shown in Fig. 4, the load at node 41 is an induction motor load while the rest of loads in the off-outage area are static loads. The parameters of the induction motor at node 41 are given in Table I. These are according to the estimated parameters of the induction motor in the physical test setup at INESC TEC.

Table I. The estimated parameters of the induction motor used in the practical test setup at INESC TEC.

| Size (W) | $R_1(\Omega)$ | $X_1(\Omega)$ | $R_2(\Omega)$ | $X_2(\Omega)$ | $X_m(\Omega)$ | $H(\sec)$ |
|---|---|---|---|---|---|---|
| 4000 | 1.44 | 2.56 | 1.37 | 2.56 | 56.17 | 0.198 |

In this paper, the original test network reported in [18] is updated adding one dispatchable DG at bus 42 with active and apparent power capacities equal to 1.08MW and 2MVA, respectively. By default, the under-voltage protection devices are installed only for the motor loads. Over-current relays exist to protect the motor loads and the root lines of each feeder. The protection devices act according to the curves given in Fig. 3. The nominal active and reactive load values for each node is according to the practical data given in [19]. In this paper, the static loads are assumed to be of constant-impedance type. Therefore, the load-voltage sensitivity coefficients are set to $K_p = 2$ and $K_q = 2$. It is assumed that each node in the network shown in Fig. 4, is equipped with a load breaker. In this test distribution network, some of the loads are considered as critical loads which have the priority factors equal to 100, whereas the priority factors of other loads are equal to 1. The critical loads are shown with '*'. in Fig. 4.

### A. Optimization Results

Three different simulation scenarios are tested in order to obtain I) the optimal status of loads in the off-outage area shown in Fig. 4 and II) the optimal active/reactive power set points of the DG at node 42 while the motor is re-energized.

Scenario I is defined assuming that the mechanical load of the induction motor has a linear torque-speed characteristic. The magnitude of the nominal mechanical torque is equal to 2 p.u. (at the synchronous speed).

The optimization results are reported in Table II. The second column and third columns provide the optimal decisions. As it can be seen, the maximum apparent power capacity of the DG (3.0 p.u.) is used. The last two columns of Table II show representative numerical information corresponding to the obtained optimal solution. This information includes the minimum margin of the current and voltage magnitudes during the motor acceleration, with respect to the over-current and under-voltage limits over all the lines and nodes equipped with protection devices.

Table II. Optimization results in case of different simulation scenarios.

| Scenarios | Load Shedding | $P^{DG}/Q^{DG}$ (p.u.) | $F^{re}$ (p.u.) | Computation time (sec) | Min. current margin (p.u.) | Min. voltage margin (p.u.) |
|---|---|---|---|---|---|---|
| I | {31,34,35,36,38,39,40} | 1.57/2.55 | 140.9 | 2.71 | 18.4 at breaker B6 | 0.0000 at node 41 |
| II | {31,33,34,36,37,38,39} | 1.8/2.4 | 200.8 | 3.20 | 18.9 at breaker B6 | 0.0000 at node 40 |
| III | {31,32,33,34,35,36, 38} | 1.89/2.32 | 210.7 | 6.125 | 19.0 at breaker B6 | 0.0013 at node 41 |

Simulation scenario II is studied in order to see how putting more voltage dip constraints would influence the optimal solution. In this regard, we assume the same simulation conditions as the ones in scenario I except that the under-voltage protections will be placed at all the nodes in the off-outage area. As shown in Table II, this change results in 69.9 p.u. increase in the reliability objective value. The computation time is still kept reasonable for the restoration problem. As expected, the minimum under-voltage margin occurs at node 40, which is the leaf node of feeder E. It is worth to highlight the fact that according to (17), it is not needed to respect the under-voltage limit at the nodes with detached loads. For example, in scenario II, although all the nodes are equipped with under-voltage protection, it is needed to respect the under-voltage constraints only at the nodes with restored loads, which are {31,33,34,36,37,38,39}.

In simulation scenario III, it is assumed that the mechanical load on the shaft of the induction motor has a fixed torque equal to 0.95 p.u. The rest of simulation conditions are the same as the ones in scenario I. As shown in Table II, although the nominal mechanical load in scenario III is so much less than the one in scenario I, more loads should be detached in order to accelerate the motor in a safe way. The reason is that under a fixed-torque mechanical load, the induction motor can accelerate only if the starting torque (electrical torque at the standstill, s=1) is larger than the mechanical torque. In order to generate this starting toque, the voltage at the motor terminal should be large enough which is obtained by detaching enough loads in the off-outage area.

The following sub-sections provide validation results. For the space limitations, the validation studies are done only for simulation scenario II. The aim is to show that the proposed semi-static optimization model represents the dynamics of the induction motor starting into the optimization problem with sufficient degree of accuracy. In this regard, the optimal decisions obtained for scenario II are deployed on the test network and then the feasibility of the induction motor starting in terms of the transient constraints is evaluated.

### B. Time-domain simulation results

As the first step of the feasibility validation, the time-domain simulation are presented. In this regard, an off-line model of the distribution network shown in Fig. 4 is built in Matlab/Simulink. In this model, the motor loads and static loads are represented by the Simulink model of the induction motor, and the impedance loads, respectively. The DG is modeled with a controllable PQ dynamic load. The optimization results of scenario II that were reported in the second and third columns of Table II are deployed on the Simulink model. Fig. 5 shows the simulation results. The results illustrate that the voltage profile at node 40 is in the safe region with respect to the under-voltage protection curve. As it can be seen, the voltage magnitude has a very narrow margin which indicates the optimality of the obtained restoration solution.

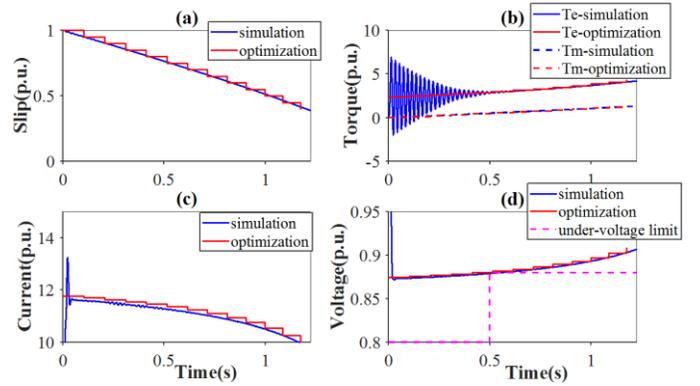

Fig. 5. Off-line simulation results for scenario II. a) slip, b) electrical and mechanical torque, c) current in breaker B6, d) voltage at node 40.

Fig. 5 also compares the electrical state profiles obtained from the optimization problem with the ones obtained from the time-domain simulation. Fig. 5.c illustrates that the initial overshoot of the motor current cannot be represented by the proposed semi-static optimization mode. As mentioned in section II, the DC term of the starting current of the induction motor is neglected in deriving the semi-static optimization model. Disregarding the very fast transients, Fig. 5 shows that the proposed semi-static model represents accurately the behavior of electrical state variables during the motor acceleration, both at the motor side and network side.

### C. Experimental test

In this section, it is aimed to check the feasibility of the optimization results using a Power Hardware-In-the-Loop (PHIL) experiment. In the test framework, instead of using an off-line simulation model of the induction motor at node 41, the actual hardware motor is integrated into a real-time simulation model of the distribution network. Therefore, the mutual interaction of the physical induction motor and the detailed model of the distribution network during the motor acceleration are considered. The test setup is implemented as shown in Fig. 7 in the smart grid laboratory at INESC TEC.

The block diagram of the laboratory test setup is depicted in Fig. 6. This PHIL test setup consists of three main parts. I) The first part includes the Matlab/Simulink network model of all the elements shown in Fig. 4 except the motor at node 41. The optimal decisions (load shedding and DG setting) obtained for simulation scenario II and given in Table II will be deployed on the simulation model of the network. The resulting simulation model is compiled, and then executed by the Real-Time

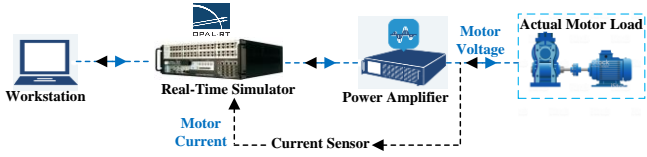

Fig. 6. Block diagram of the PHIL test setup.

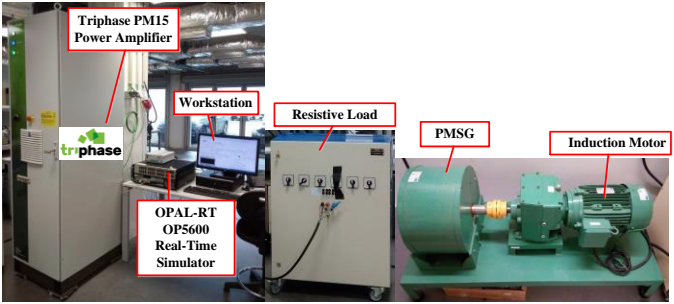

Fig. 7. The experimental PHIL test setup in the smart grid laboratory at INESC TEC.

Simulator with the aim of emulating dynamics of the simulated system in real-time. Once the compiled model is executed by the Real-Time Simulator, the obtained voltage signal at the bus connected to the motor load is exported to the power amplifier in the form of a digital signal. The Real-Time Simulator used in this test setup is an OPAL-RT OP5600 shown in Fig. 6. II) The second part includes the Power Amplifier which magnifies the motor voltage signal received from the Real-Time Simulator. This real voltage is applied to the physical motor hardware. Then the current measured at the output of the Power Amplifier (the motor current) is sent back to the Real-Time Simulator. Therefore, the loop of the electric circuit encompassing the hardware and software devices is closed. As shown in Fig. 7, the Power Amplifier used in this test setup is a TriPhase PM15. This Power Amplifier has a nominal voltage equals to 400V and it can tolerate up to 30A (peak currents per phase). The communication between the Real-Time Simulator and the Power Amplifier is established using a duplex fiber optic connection. III) The third part of the PHIL test setup is the physical hardware. The hardware system in this experiment includes a mechanical switch and an induction motor that needs to be accelerated under a mechanical load. As shown in Fig. 7, in order to emulate a mechanical load with a linear torque-speed characteristic, a Permanent Magnet Synchronous Generator (PMSG) is used connected to a resistive load. The torque that is imposed on the shaft by this PMSG changes almost linearly with rotational speed. The parameters of the induction motor shown in Fig. 7 are as given in Table I.

In order to execute the PHIL experiment, first, the closed-loop stability of the PHIL configuration have to be ensured. Then, the PHIL is initiated while the hardware switch is open and the hardware motor is de-energized. At the end, the hardware switch will be closed in order to accelerate the induction motor. Fig. 8 validates the feasibility of the optimal solution found in scenario II regarding the under-voltage and over-current constraints.

## V. CONCLUSION

This paper proposed a semi-static analytical optimization model for the starting of the induction motors. This model was used to develop a mixed-integer convex formulation for the restoration problem. In this restoration problem, it is aimed to re-energize the induction motor following a fault in a safe way by making optimal decisions for the I) load pickup sequence and II) DG power set points. The overall optimization problem was formulated in terms of a MISOCP problem. This optimization problem was solved using the Gurobi solver for a test distribution network. It was shown that the optimal solution is obtained within a short time for different simulation scenarios. The very narrow margin obtained for the voltage magnitude with respect to the under-voltage limit indicates the high quality of the restoration solution. The feasibility of the solution was verified using off-line time-domain simulations and also using a Power-Hardware-In-the-Loop test experiment. The comparison of the state variables obtained from the optimization problem with the ones obtained from the simulation studies illustrates the high accuracy of the proposed semi-static optimization problem.

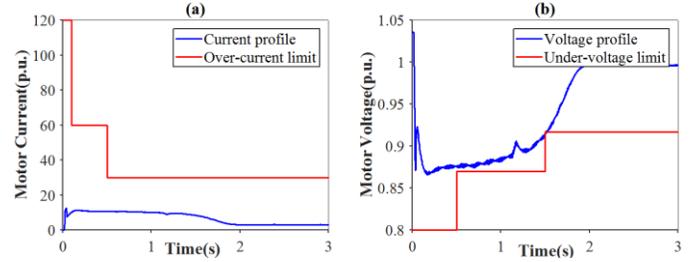

Fig. 8. The results of the PHIL test experiment in scenario II, a) starting current magnitude, b) terminal voltage magnitude.

## I. ACKNOWLEDGMENT

The authors gratefully acknowledge financial support of the Qatar Environment and Energy Research Institute (QEERI).

## III. Appendices

### A. Slip discretization

In this section, it is aimed to justify the assumption made in section II that the slip is fixed at each step. Equation (19) is the first-order approximation of the tailor expression for the motor slip as a function of time. This approximation is made for the time interval $[t_k, t_k + \Delta t_k]$, where $t_k$ and $\Delta t_k$ are the starting time and the time length of step $k$, respectively.

$$S(t) \approx S_k + \frac{dS}{dt}(t - t_k) \quad \forall t \in [t_k, t_k + \Delta t_k] \quad (19)$$

The time derivative of slip can be expressed according to the dynamic motion equation (20).

$$\frac{dS}{dt} = \frac{1}{2H}\left(T^{ele} - T^{mec} - K_D(1 - S)\right) \quad (20)$$

As it can be understood from (19), the time derivative of slip should be small enough in order to assume that the slip equals to $S_k$ during each step $k$. According to (20), it will be the case when the total inertia of the induction motor is large. This condition holds when we study the starting of large induction motors under mechanical loads with high moment of inertia or when we study an aggregated model of multiple induction motors. According to (19), for modeling low-inertia induction motors, we have to assume more number of slip steps so that the time length of each step $k$ ($\Delta t_k$) will be small enough. Therefore, we can assume that the second expression in the right hand side of (19) is almost zero and therefore the slip can be approximated to $S_k$ during the step $k$.

### B. Piecewise Linear Approximation

Let $f(x)$ is an arbitrary continuous function with the domain of $[x_1, x_n]$. Fig. 9 shows $\tilde{f}(x)$, in red, as the piece-wise linear approximation of function $f(x)$. Let divide the domain of function $f$ by $n$ break-points $x_1, x_2, \dots, x_n$. A break-point is defined as the point where the slop of $\tilde{f}(x)$ changes. Constraints (21)-(24) formulate $\tilde{f}(x)$ by introducing $n$ non-negative auxiliary variables $\lambda_i$. Any point in between two breakpoints is a weighted sum of these two breakpoints. Therefore, as expressed in (21), a given $x$ value that is within the domain of function $f$ can be expressed as a weighted sum of variables $\lambda_i$. Using the same auxiliary variables, function $\tilde{f}(x)$ is formulated as in (22). Therefore, $\tilde{f}(x)$ consists of straight lines between successive breakpoints. The sum of all auxiliary variables $\lambda_i$ should be one (23). Moreover, out of all $\lambda_i$ variables, at most two successive variables can be non-zero (24). This latter constraint is in the form of a common restriction that is known as Special Ordered Set-2 (SOS2) constraint [5]. Most of the commercial solvers have the feature to account for these types of constraints during the Branch-and-Bound search algorithm. The more number of break points, the more accurate approximation will be obtained. However, the more number of breakpoints will require more number of variables which increases the computation burden.

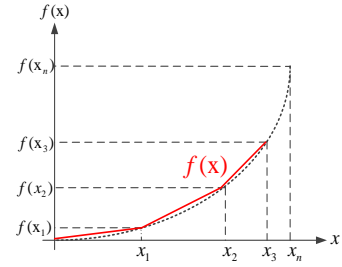

Fig. 9. Piece-wise linear approximation of an arbitrary continuous function $f(x)$

$$x = \sum_{i=1}^{n} \lambda_i x_i \quad (21)$$

$$\tilde{f}(x) = \sum_{i=1}^{n} \lambda_i f(x_i) \quad (22)$$

$$\sum_{i=1}^{n} \lambda_i = 1 \quad (23)$$

$$\lambda_i \geq 0 \; : SOS2 \quad (24)$$

### C. Elimination of product of variables

In this section, a method is provided for the linearization of the constraints which incorporate a product of two variables. The product of two variables $x_1$ and $x_2$ can be replaced by one new variable y, on which the constraints given in (25) and (26) are imposed. It is assumed that $x_1$ is a binary variable and $x_2$ is a positive continuous variable, for which $0 \leq x_2 \leq u$ holds. If $x_1$ is zero, $y$ must be zero according to (25), otherwise if $x_1$ is one, (26) forces $y$ to take the value of $x_2$.

$$0 \leq y \leq u x_1 \quad (25)$$

$$x_2 - u(1 - x_1) \leq y \leq x_2 \quad (26)$$